\documentclass{article}

\newtheorem{theorem}{Theorem}

\newtheorem{lemma}{Lemma}

\newtheorem{corollary}{Corollary}

\begin{document}

\bibliographystyle{amsplain}

\title{On the existence of shortest directed networks}
\author{Konrad J. Swanepoel\\ 
	Department of Mathematics and Applied Mathematics \\ 
	University of Pretoria,
	Pretoria 0002 \\ South Africa \\
	E-mail: \texttt{konrad@math.up.ac.za}}
\date{}
\maketitle

\begin{abstract}
A \emph{directed network} connecting a set $A$ to a set $B$ is a digraph containing an $a$-$b$ path for each $a\in A$ and $b\in B$.
Vertices in the directed network not in $A\cup B$ are \emph{Steiner points}.
We show that in a finitely compact metric space in which geodesics exist, any two finite sets $A$ and $B$ are connected by a shortest directed network.
We also bound the number of Steiner points by a function of the sizes of $A$ and $B$.
Previously, such an existence result was known only for the Euclidean plane [M. Alfaro, Pacific J. Math.\ \textbf{167} (1995) 201--214].
The main difficulty is that, unlike the undirected case (Steiner minimal trees), the underlying graphs need not be acyclic.

Existence in the undirected case was first shown by E. J. Cockayne [Canad.\ Math.\ Bull.\ \textbf{10} (1967) 431--450].

\end{abstract}

\quote{\emph{Dedicated to Ernie Cockayne on the occasion of his 60th birthday}}

\section{Introduction}
Let $(X,\rho)$ be a metric space, which we assume to be \emph{finitely compact}, i.e.\ closed and bounded sets are compact, and in which we furthermore assume that any two points $x,y\in X$ are connected by a \emph{geodesic}, i.e.\ an arc of length $\rho(x,y)$.

In this general context, Cockayne \cite{Cockayne} was the first to show that for any finite set $A\subset X$, there exists a \emph{shortest} undirected network (i.e.\ a Steiner minimal tree) connecting the points of $A$.
In this paper we consider the case of directed networks (defined below).
Unlike the directed case, there can be cycles in the network (see \cite{Alfaro} for examples), which presents a problem in bounding the number of Steiner points.
We give an upper bound for the number of Steiner points depending only on the size of the set of points being connected, and from this deduce our existence result.
Previously, such an existence result was shown only in the Euclidean plane \cite{Alfaro}, with the proof depending heavily on facts of plane geometry, making it unsuitable for generalisation even to Euclidean space.
Our approach is entirely combinatorial.

In the next section we state our terminology and formulate our results.
In Section~\ref{proofsection} we prove the main theorem, and in Section~\ref{remarksection} we discuss related problems and remark on computational aspects.

\section{Terminology and Results}
Our digraph terminology follows \cite{CL}.
Let $A$ and $B$ be finite subsets of the metric space $(X,\rho)$.
A \emph{directed network} connecting $A$ to $B$, or \emph{$(A,B)$-network} for short, is a digraph $G=(V,E)$ such that $A\cup B\subseteq V\subseteq X$, and there is a directed $a$-$b$ path for each $a\in A$ and $b\in B$.
The \emph{length} of $G$ is
\[ \ell(G)=\sum_{(x,y)\in E} \rho(x,y).\]
We call the vertices in $A$ \emph{sources}, the vertices in $B$ \emph{sinks}, and (following tradition) the vertices in $V\setminus(A\cup B)$ \emph{Steiner points}.
A \emph{shortest} $(A,B)$-network in $X$ is an $(A,B)$-network of minimum length, provided it exists.

If a Steiner point in an $(A,B)$-network has at most two neighbours, it and its incident edges may be removed and (possibly) replaced by edges between the neighbours, to obtain a new $(A,B)$-network not longer than the original network.
This is easily verified by considering the various cases and using the triangle inequality.
A \emph{simple} $(A,B)$-network is an $(A,B)$-network in which each Steiner point has at least three neighbours.
Clearly, if there exists a shortest $(A,B)$-network in $X$, there also exists a simple shortest one.

\begin{theorem}
In any metric space $(X,\rho)$, if $A,B\subset X$, then any simple shortest $(A,B)$-network has at most $O(m^2n+mn^2)$ Steiner points, where $m=|A|$ and $n=|B|$.
\end{theorem}

The proof is in the next section.
To derive an existence result from this theorem, we need a compactness argument, encapsulated as follows:

\begin{lemma}
Let $(X,\rho)$ be a finitely compact metric space in which any two points are connected by a geodesic.
Let $A,B\subset X$ be finite sets and $s$ a positive integer.
Then, among all $(A,B)$-networks with at most $s$ Steiner points, there is a shortest one.
\end{lemma}

\begin{corollary}
Let $(X,\rho)$ be a finitely compact metric space in which any two points are connected by a geodesic.
Then, for any finite $A,B\subset X$, there exists a shortest $(A,B)$-network.
\end{corollary}

\textbf{Proof of Lemma} \,\,
As this is a standard compactness argument, we only sketch the proof.
Take a sequence of $(A,B)$-networks $G_i$ such that
\[ \lim_{i\to\infty} \ell(G_i) = \inf \ell(G),\]
where the infimum is taken over all $(A,B)$-networks $G$ with at most $s$ Steiner points.
The $G_i$'s are all contained in some bounded subset of $X$, so we may assume without loss that $X$ is compact.
We may take a subsequence such that the abstract digraph structure of all $G_i$'s are the same (say $G$), since there are only finitely many digraphs with at most $|A|+|B|+s$ vertices.
We may again take subsequences until all Steiner points converge, since $X$ is compact.
In the limit we obtain an $(A,B)$-network (with underlying digraph a contraction of $G$) of length $\lim_i \ell(G_i)$.

\section{Bounding the number of Steiner points}\label{proofsection}
Given a set of sources $A$ and a set of sinks $B$ of sizes $|A|=m$ and $|B|=n$ in the metric space, we let $G$ be any simple shortest $(A,B)$-network.
Note that by minimality, $G$ is covered by all $a$-$b$ paths where $a\in A,\,b\in B$.
We now show that each $a$-$b$ path contains at most $O(m+n)$ points, thereby proving the Theorem.

It is sufficient to prove that the longest $a$-$b$ path has at most $O(m+n)$ points, where we take the longest path over all $a$-$b$ paths of all simple shortest $(A,B)$-networks.
Let $P=a_0x_1\dots x_kb_0$ be such a longest $a$-$b$ path.
We may assume that $G$ is the simple shortest $(A,B)$-network containing this path.
Also, let $x_0=a_0$ and $x_{k+1}=b_0$.

For each source $a$, let $x(a)$ be the first point in $P$ such that there exists an $a$-$x(a)$ path in $G$, and fix such a path $P(a)$.
Such an $x(a)$ always exists, since there is at least an $a$-$b_0$ path.
Similarly, let $y(b)$ be the last point in $P$ such that there exists a $y(b)$-$b$ path $Q(b)$ in $G$, for each sink $b$.

For each source $a$ and sink $b$, fix an $a$-$b$ path $P(a,b)$.
If $P(a,b)$ contains an $x_i$ and an $x_j$ $(i<j)$ with $x_i$ appearing before $x_j$ in $P(a,b)$, we may replace the $x_i$-$x_j$ subpath of $P(a,b)$ by the subpath $x_ix_{i+1}\dots x_j$ of $P$.

Also, we may replace the initial segment of $P(a,b)$ from $a$ to the first point of $P$ on $P(a,b)$ by $P(a)$, as well as the final segment of $P(a,b)$ from the last point of $P$ on $P(a,b)$ to $b$ by $Q(b)$.

We may therefore assume that $P(a,b)$ is either disjoint with $P$, or 
\begin{itemize}
\item starts off with $P(a)$, 
\item then consists of subpaths of $P$ and $x_i$-$x_j$ paths $(i>j)$ edge-disjoint with $P$, which we call \emph{$(i,j)$-jumps}, 
\item and then ends with $Q(b)$.
\end{itemize}
In particular, $P(a_0,b_0)=P$.

Let $J$ be the set of all $(i,j)$-jumps appearing in all $P(a,b)$'s.
By minimality of $G$, $G$ is covered by the union of all $(i,j)$-jumps, all $P(a)$'s, all $Q(b)$'s, and $P$.
Let $I$ be a minimal subset of $J$ such that $J$ may be replaced by $I$ in the above union, and $G$ is still covered.
We modify each $P(a,b)$ so as to use only $(i,j)$-jumps from $I$, by replacing each $(i,j)$-jump in $J\setminus I$ by an $x_i$-$x_j$ path consisting of $P(a)$'s, $Q(b)$'s, subpaths of $P$ and $(i,j)$-jumps from $I$ (such paths existing because $G$ is still covered).

We now show that each vertex of $P$ is either some $x(a)$, some $y(b)$ or an endpoint of some $(i,j)$-jump from $I$.
Consider any $x_t$ $(1\leq t\leq k)$.
Since $G$ is simple, $x_t$ is either a source or a sink, hence an $x(a)$ or $y(b)$, or a Steiner point, in which case it is connected to some vertex $c\neq x_{t-1},x_{t+1}$.
Since $x_tc$ is not a redundant edge, it must be contained in all paths connecting some source $a$ to some sink $b$, hence must be in $P(a,b)$.
Thus $c$ is contained in some $(i,j)$-jump of $P(a,b)$.
Since an $(i,j)$-jump is edge-disjoint with $P$, $c$ must be an endpoint.

Thus the number of vertices of $P$ is bounded above by $m+n+2|I|$.
We now bound $|I|$ from above.
Note that by minimality, for any $(i_1,j_1)$-jump and $(i_2,j_2)$-jump in $I$ we have $i_1\neq i_2$ and $j_1\neq j_2$.
Denote the unique $(i,j)$-jump by $(i,j)$.
Define a relation $\succ$ between $I$ and
\[X=\{x(a):a\in A\}\cup\{y(b):b\in B\}\]
by
\[(i,j) \succ x_t \mbox{ iff } j\leq t\leq i.\]
Call $(i_1,j_1)$ and $(i_2,j_2)$ \emph{consecutive} if there is no $(i_3,j_3)\in I$ with $i_3$ between $i_1$ and $i_2$.
The relation $\succ$ has the following two properties:
\begin{enumerate}
\item For each $x_t\in X$ there are at most two $(i,j)\in I$ such that $(i,j)\succ x_t$.
\item For any two consecutive $(i,j),(k,\ell)\in I$ there is at least one $x_t\in X$ such that $(i,j)\succ x_t$ or $(k,\ell)\succ x_t$.
\end{enumerate}

From these two properties it follows that $|I|\leq 4(m+n)+1$, hence the number of vertices of $P$ is at most $9(m+n)+2$.

It remains to verify the above two properties.

For the first, suppose that $(i_1,j_1),(i_2,j_2),(i_3,j_3)\succ x_t$, with $i_1> i_2> i_3$.
If $j_1< j_2$, then $(i_2,j_2)$ is redundant.
Therefore, $j_1>j_2$, and similarly, $j_2>j_3$. But then $(i_2,j_2)$ is again redundant.

For the second property, note that by minimality of $I$, if $(i,j)$ and $(k,\ell)$ are consecutive with $i>k$, then $j>\ell$.
There are now two cases to consider:

$j<k:$
If $k-j=1$, then we may change the directions of the $(k,\ell)$-jump and of $x_\ell\dots x_j$ and discard $x_jx_k$ to obtain a shorter $(A,B)$-network.
Therefore, $k-j\geq 2$.
If $x_{j+1}$ is an endpoint of some $(i',j')$, then we obtain that $I$ is not minimal.
Thus $(i,j)\succ x_{j+1}\in X$.

$j\geq k:$
If $i-j\geq 2$ and $x_{i-1}\not\in X$, then $x_{i-1}$ is an endpoint of an $(i',i-1)$ with $i'>i$.
The previous case then provides a contradiction.
Thus $(i,j)\succ x_{i-1}\in X$.
Similarly, if $k-\ell\geq 2$, then $(k,\ell)\succ x_{\ell+1}\in X$.
So we may assume $i-j=1$ and $k-\ell=1$.
If $(i,j)$ has more than one edge, we may redirect $G$ to give an $(A,B)$-network with a longer $a_0$-$b_0$ path, a contradiction.
Therefore, $(i,j)$ has only one edge.
Similarly, $(k,\ell)$ has only one edge.
Since $G$ is simple, $x_j\in X$.

This proves the Theorem.

\section{Concluding remarks}\label{remarksection}
We have made no attempt at finding the best order for the number of Steiner points. 
For example, with some more effort the constant $9$ in the above upper bound for the number of points on an $a$-$b$ path can be lowered.
We also mention that there are metric spaces with sets needing at least $m+n$ Steiner points.
We believe that the best upper bound should be less than cubic. 
Finding a better upper bound will be crucial in designing an algorithm that will be effective for at least small sets of points.

As expected, finding a shortest $(A,B)$-network in a digraph is NP-complete.
The Minimal Equivalent Digraph (MED) problem is that of finding a minimal subgraph of a graph, such that if any two vertices are connected by a path in the original digraph, they are still connected in the subdigraph.
This problem is known to be NP-complete \cite{Sahni}, and it is easily seen that the MED problem is polynomial time reducible to finding a shortest $(A,B)$-network in a digraph.

A related problem in Operations Research is the point-to-point connection problem (see e.g.\ \cite{NF}), where $A=\{a_1,a_2,\dots,a_n\}, B=\{b_1,b_2,\dots,b_n\}$, and we require only that there are $a_i$-$b_i$ paths for each $i=1,\dots,n$.
The proof of the Theorem shows that in this case a simple network has at most $O(n^2)$ Steiner points.

\providecommand{\bysame}{\leavevmode\hbox to3em{\hrulefill}\thinspace}

\end{document}